\documentclass[12pt,leqno,draft]{article}

\newtheorem{theorem}{Theorem}
\newtheorem{lemma}[theorem]{Lemma}
\newtheorem{proposition}[theorem]{Proposition}
\newtheorem{definition}[theorem]{Definition}
\newtheorem{corollary}[theorem]{Corollary}

\newcommand{\begintheorem}{\addtocounter{equation}{1}\begin{theorem}}
\newcommand{\beginlemma}{\addtocounter{equation}{1}\begin{lemma}}
\newcommand{\beginproposition}{\addtocounter{equation}{1}\begin{proposition}}
\newcommand{\begindefinition}{\addtocounter{equation}{1}\begin{definition}}
\newcommand{\begincorollary}{\addtocounter{equation}{1}\begin{corollary}}

\begin{document}

\title{Some remarks about curves in metric spaces}

\author{Stephen William Semmes	\\
	Rice University		\\
	Houston, Texas}

\date{}

\maketitle

	Let $(M, d(x, y))$ be a metric space.  Thus $M$ is a nonemtpy
set and $d(x, y)$ is a nonnegative real-valued function defined for
$x, y \in M$ such that $d(x, y) = 0$ if and only if $x = y$,
$d(x, y) = d(y, x)$ for all $x, y \in M$, and
\begin{equation}
	d(x, z) \le d(x, y) + d(y, z)
\end{equation}
for all $x, y, z \in M$.  Of course this last condition is known as
the triangle inequality.

	We can weaken the triangle inequality to the requirement
that there is a positive real number $C$ such that
\begin{equation}
	d(x, z) \le C \, (d(x, y) + d(y, z))
\end{equation}
for all $x, y, z \in M$.  In this case, with the other conditions
as before, we say that $d(x, y)$ is a quasimetric on $M$.  A stronger
version of the triangle inequality asks that
\begin{equation}
	d(x, z) \le \max (d(x, y), d(y, z)),
\end{equation}
and when this happens we say that $d(x, y)$ is an ultrametric.
One can check that an ultrametric space is totally disconnected,
which is to say that it does not contain a connected subset
with more than two elements.

	Let us say that a subset $E$ of a metric space $(M, d(x, y))$
is chain connected if for every pair of points $u, v \in E$ and every
$\epsilon > 0$ there is a finite chain $w_1, \ldots, w_l$ of points in
$E$ such that $w_1 = u$, $w_l = v$, and $d(w_j, w_{j+1}) < \epsilon$
for all $1 \le j < l$.  For any subset $E$ of $M$ and $\epsilon > 0$,
if $u$ is an element of $E$ and $E_1(u)$ is the set of points in $E$
which can be connected to $u$ by a finite $\epsilon$-chain of points
in $E$ of this type, and if $E_2(u)$ consists of the remaining points
in $E$, then the distance between every element of $E_1(u)$ and
$E_2(u)$ is at least $\epsilon$.  As a consequence, if $E$ is connected
in the usual sense, then $E$ is chain connected.  One can check that
the converse holds when $E$ is compact.

	By a path in a metric space we mean a continuous mapping from
a closed and bounded interval $[a, b]$ in the real line into $M$.  A
subset $E$ of $M$ is said to be pathwise connected if for every pair
of points $u, v \in E$ there is a continuous path contained in $E$
which begins at $u$ and ends at $v$.  A pathwise-connected set is
connected, but the converse does not work in general, even for compact
subsets of ${\bf R}^2$.  As in the previous paragraph, connectedness
implies chain connectedness, which is somewhat like path connectedness,
but without much information on the complexity of the chains.

	One can consider more refined notions of chain connectedness
and pathwise connectedness with controls on the complexity of the
chains or paths.  For that matter one can view $\epsilon$-chains as a
kind of generalization of paths, defined on a discrete set of points
in the real line.  For instance one might choose the points in the
domain so that their incremental distances are the same as the
corresponding points in the metric space.

	The types of controls that one might consider for chains or
paths are closely related to the kind of metric being used.  If $(M,
d(x, y))$ is a metric space and $a$ is a positive real number, one can
define a new distance function $\rho(x, y)$ on $M$ by
\begin{equation}
	\rho(x, y) = d(x, y)^a.
\end{equation}
If $0 < a < 1$, one can check that this defines a metric on $M$,
which we may call the snowflake transform of order $a$ of $d(x, y)$.
If $a > 1$, then $\rho(x, y)$ is still a quasimetric on $M$.
If $d(x, y)$ happens to be an ultrametric, then $\rho(x, y)$
is also an ultrametric for all $a > 0$.

	Let us mention a very nice converse result from \cite{2}.
Namely, if $\rho(x, y)$ is a quasimetric on $M$, then there is a
metric $\delta(x, y)$ on $M$ and positive real numbers $\eta$, $C$
such that $C^{-1} \, \delta(x, y) \le \rho(x, y)^\eta \le C \,
\delta(x, y)$ for all $x, y \in M$.  Thus quasimetrics can always be
approximated by ordinary metrics in this manner.

	Suppose that $(M_1, d_1(x, y))$ and $(M_2, d_2(u, v))$ are
metric spaces, or even quasimetric spaces, and let $f$ be a mapping
from $M_1$ to $M_2$.  We say that $f$ is Lipschitz of order $a$
for some positive real number $a$ is there is a positive real number
$L$ such that
\begin{equation}
	d_2(f(x), f(y)) \le L \, d_1(x, y)^a
\end{equation}
for all $x, y \in M$.  This parameter $a$ is closely related to the
exponents of distance functions discussed earlier, because one can
change $a$ automatically by replacing $d_1(x, y)$ or $d_2(u, v)$
by positive powers of themselves.

	An important feature of metric spaces is that they always have
a rich supply of real-valued Lipschitz functions of order $1$.
To be more precise, if $(M, d(x, y))$ is a metric space and
$p$ is any element of $M$, then the function $f_p(x) = d(x, p)$
is Lipschitz of order $1$, with constant $L = 1$.  This can be
verified using the triangle inequality, and it does not work in
general for quasimetrics.  Of course we use the standard metric
on the real line for the range of these functions.

	We can use Lipschitz conditions to control the complexity of
curves in metric spaces, or also chains of points by viewing them as
mappings from discrete subsets of the real line into the metric space.
If $p(t)$ is a Lipschitz mapping of order $1$ from the unit interval
$[0, 1]$ in the real line into a metric space $(M, d(x, y))$, then it
is reasonable to say that the path has finite length less than or
equal to the Lipschitz constant of the mapping.  In general it may be
possible to connect a pair of points in a subset $E$ of a metric space
$M$ by a continuous path in $E$, and one which is even Lipschitz of
some orders, and not Lipschitz of other orders.

	If $(M, d(x, y))$ is a metric space and $a$ is a real number
such that $a > 1$, then any continuous mapping from an interval
in the real line into $M$ which is Lipschitz of order $a$ is constant.
When $M$ is the real line, with the usual metric, this follows
from the observation that a real-valued function on an intevral
which is Lipschitz of order strictly larger then $1$ has derivative
$0$ everywhere.  In general one can reduce to this case by mapping
the curve from $M$ into the real line using a real-valued Lipschitz
function.

	Now suppose that $(M, d(x, y))$ is a metric space and $b$ is a
positive real number with $b < 1$, and consider the snowflake metric
$\rho(x, y) = d(x, y)^b$.  Any Lipschitz mapping of order $1$ from an
interval in the real line into $(M, \rho(x, y))$ is the same as a
Lipschitz mapping of order $a = 1/b > 1$ into $(M, d(x,y))$, and hence
is constant.  There may be curves defined by Lipschitz mappings of
order $\le b$, depending on the gometry of $M$.

	There are a lot of classical topics in geometric topology
related to dimensions and embeddings, as in \cite{1}.  In particular
let us mention the famous examples of the Sierpinski gasket and carpet
and the Menger sponge.  The first two are compact subsets of ${\bf
R}^2$ while the third is a compact subset of ${\bf R}^3$.  Each has
topological dimension $1$, and is also pathwise connected.

	In fact these well-known fractal sets also have a lot of
nice curves of finite length.  From a purely topological point of
view this might be considered as an extra bonus.  They do not have
any snowflaking, and they do not need any.  Of course there are also
matters of self-similarity, nice measures on them, etc.

	Purely topological aspects of spaces like these and related
mappings have been studied quite a bit.  One might also mention other
kinds of compact connected sets such as Bing's pseudo-arc, which is a
lot like a continuous arc but is not an arc, and has other special
features too.  There are a lot of tricky properties of spaces like
these, along the lines of what might be mapped where satisfying
such-and-such conditions.  Additional restrictions on complexity
such as those given by Lipschitz classes lead to a lot of new
questions.

\end{document}